\documentclass[12pt]{amsart}

\def\@typesizes{
       \or{5}{6.5}\or{6}{7.5}\or{7}{8.5}\or{8}{11}\or{9}{12}
       \or{10}{13}
       \or{\@xipt}{14}\or{\@xiipt}{15}\or{\@xivpt}{18}
       \or{\@xviipt}{20}\or{\@xxpt}{24}}

\usepackage{amsthm}
\usepackage{amsmath}
\usepackage{amssymb}
\usepackage{graphicx}
\usepackage{enumerate}
\usepackage{xcolor}
\usepackage{tikz}

\allowdisplaybreaks

\numberwithin{equation}{section}
\numberwithin{figure}{section}

\theoremstyle{plain}
\newtheorem{theorem}{ Theorem}[section]
\newtheorem{proposition}[theorem]{ Proposition}
\newtheorem{lemma}[theorem]{ Lemma}
\newtheorem{corollary}[theorem]{ Corollary}
\newtheorem{example}[theorem]{ Example}
\newtheorem{remark}[theorem]{ Remark}
\newtheorem{definition}[theorem]{ Definition}
\newtheorem{conjecture}{ Conjecture}

\def\BET{\begin{theorem}}
\def\ENT{\end{theorem}}
\def\BEP{\begin{proposition}}
\def\ENP{\end{proposition}}
\def\BEL{\begin{lemma}}
\def\ENL{\end{lemma}}
\def\BEC{\begin{corollary}}
\def\ENC{\end{corollary}}
\def\BEE{\begin{example} \rm}
\def\ENE{\end{example}}
\def\BER{\begin{remark} \rm}
\def\ENR{\end{remark}}
\def\BED{\begin{definition} \rm}
\def\END{\end{definition}}
\def\BECJ{\begin{conjecture}}
\def\ENCJ{\end{conjecture}}

\def\bea{\begin{eqnarray}}
\def\eea{\end{eqnarray}}

\def\beas{\begin{eqnarray*}}
\def\eeas{\end{eqnarray*}}

\def\beq{\begin{equation}}
\def\eeq{\end{equation}}

\def\beal{\begin{align*}}

\def\eeal{ \end{align*} }

\def\roweq{\nonumber \\ &=& }
\def\rowleq{\nonumber \\  & \leq & }
\def\rowgeq{\nonumber \\ & \geq & }

\def\rowpl{\nonumber \\ & &  +  }

\def\bbB{{\mathbb B}}
\def\bbC{{\mathbb C}}
\def\bbD{{\mathbb D}}

\def\bbN{{\mathbb N}}

\def\bbR{{\mathbb R}}

\def\cA{{\mathcal A}}

\DeclareMathOperator{\wind}{wind}
\DeclareMathOperator{\ind}{ind}
\DeclareMathOperator{\BMO}{BMO}
\DeclareMathOperator{\VMO}{VMO}
\DeclareMathOperator{\BWMO}{BWMO}
\DeclareMathOperator{\VWMO}{VWMO}
\DeclareMathOperator{\BO}{BO}
\DeclareMathOperator{\VO}{VO}
\DeclareMathOperator{\BA}{BA}
\DeclareMathOperator{\VA}{VA}

\thanks{This project has received funding from the European Union’s Horizon 2020 research and innovation programme under the Marie Sklodowska-Curie grant agreement No 844451. The second author was supported in part by Engineering and Physical Sciences Research Council grant EP/T008636/1.}
\begin{document}

\title[Weak BMO and Toeplitz operators]{Weak BMO and Toeplitz operators\\on Bergman spaces}

\author{Jari Taskinen}
\address{Department of Mathematics and Statistics, University of Helsinki, Finland}
\email{jari.taskinen@helsinki.fi}

\author{Jani A. Virtanen}
\address{Department of Mathematics and Statistics, University of Reading, England}
\email{j.a.virtanen@reading.ac.uk}



\begin{abstract}
Inspired by our previous work on the boundedness of Toeplitz operators, 
we introduce weak $\BMO$ and $\VMO$ type conditions, denoted by
BWMO and VWMO, respectively, for functions on the open unit disc of the 
complex plane. We show that the average function of 
a function $f \in$ BWMO is boundedly oscillating, and the analogous
result holds for $f \in$ VWMO. The result is applied for generalizations
of known results on the essential spectra and norms of Toeplitz 
operators. Finally, we provide examples of functions satisfying 
the $\VWMO$ condition which are not in the classical $\VMO$ or even in 
$\BMO$.
\end{abstract}
\maketitle

\section{Introduction and main results}
\label{sec1}

Consider the Banach space $L^p := ( L^p(\bbD, dA), \Vert \cdot \Vert_p)$,
where $1< p < \infty$ and $dA$ is the normalized area measure
on the unit disc $\bbD$ of the complex plane $\bbC$, and the \emph{Bergman space}
$A^p$, which is the closed subspace of $L^p$ consisting of analytic functions. 
The \emph{Bergman projection} $P$ is the orthogonal projection of
$L^2$ onto $A^2$, and it  has the integral representation
\beas
Pg(z) = \int\limits_\bbD K_\zeta(z) g(\zeta)  dA (\zeta) \ , \ \ \ \mbox{where}
\ K_\zeta(z) = \frac{1 }{ ( 1 - z \bar \zeta)^2 }, \ z ,\zeta \in \bbD,  
\eeas
is the \emph{Bergman kernel}.
It is also known to be a bounded projection of $L^p$ onto $A^p$ for every 
$1<p<\infty$.  For an integrable function $f : \bbD \to \bbC$ and, say, 
bounded analytic functions $g$,  the \emph{Toeplitz operator}
$T_f$ with \emph{symbol} $f$ is defined by 
\bea
T_f g = P(fg) = \int\limits_\bbD \frac{f(\zeta)g(\zeta)  }{ ( 1 - z \bar \zeta)^2 } dA (\zeta).
\label{0.0}
\eea
A related class of operators consists of {\it Hankel operators} $H_f : 
A^p \to L^p$ defined by
$$
	H_fg  = (I-P)(fg),
$$
where $I$ is the identity operator and $I-P$ is the complementary 
projection of $P$. Notice that the boundedness of $P$ implies that both $T_f : A^p \to A^p$ and $H_f : A^p \to L^p$ are bounded whenever the symbol $f$ is in $L^\infty$.

In \cite{TV} and \cite{TV2}, we studied a generalized definition of 
Toeplitz operators with locally integrable symbols $f$ 
satisfying a weak ``averaging'' condition (see \eqref{4} below), and 
showed that one 
can define $T_f = \lim_{\rho \to 1} T_{f_\rho}$, where $f_\rho = 
\chi_\rho f$ and $\chi_\rho$ is the characteristic function of the 
compact set $\{ z \in \bbD \, : \, |z| \leq \rho\}$, $\rho < 1$. 
The limit converges in the strong operator topology and the generalized 
definition coincides with \eqref{0.0}, whenever the latter makes sense.
It was recently
shown in \cite{YZ} that the condition is however not necessary for the boundedness of $T_f$. 

Here, our aim is to apply the same idea 
to introduce new \emph{weak} $\BMO$ and $\VMO$ type conditions $\BWMO$ 
and $\VWMO$: we replace the standard definition  of  $\BMO$ (and $\VMO$) 
by the above described  weak averaging condition. It is quite clear
there are functions which belong  to $\VWMO$ but are not in $\VMO$ or not even 
in $\BMO$. We will exhibit concrete examples in Section \ref{sec4}, see 
Example~\ref{ex5}. Our new definition leads to the 
following  results. First, we prove that whenever $f$ satisfies the  
$\BWMO$ (or $\VWMO$) condition, then the average function $\widehat f$ 
belongs to the space $\BO$ (or in $\VO$) of functions of 
bounded (or vanishing) oscillation---see Theorem~\ref{prop10}. This allows us to extend the standard results on the essential spectra and Fredholm properties of Toeplitz operators $T_f$ from the case $f\in L^\infty\cap \VMO$ (see~\cite{Hagger1, PeVi} and the references therein) to integrable, not necessarily bounded, symbols in $\VWMO$.

The weak conditions are broader in scope, and should have more applications, which we hope to demonstrate in future work on Toeplitz and related operators.

\section{Preliminaries}
\label{sec2}

In this section, we explain the notions used in the paper and recall a 
number of basic results that we need in the subsequent sections. 

The space of bounded mean oscillation $\BMO^p$ provides a class of  
symbols $f$, strictly larger than $L^\infty$, for which bounded (or 
compact) Toeplitz operators can be characterized in terms of the boundary 
behavior of the Berezin transform $\widetilde f$. Similarly, its closed 
subspace of vanishing mean oscillation $\VMO^p$ plays an important role 
in the study of other (spectral) properties of Toeplitz operators. 
Let $r>0$, $1 \leq p < \infty$ and $f$ be a locally $L^p$-integrable function on $\bbD$. We 
say that $f$ is of \emph{bounded mean oscillation}, and write $f\in 
\BMO^p_r$, if
\begin{equation}\label{2.1}
\sup_{z\in \bbD} \frac1{|D(z,r)|} \int\limits_{D(z,r)} 
|f(\zeta) - \widehat f_r(z)|^p\,dA(\zeta)  < \infty,
\end{equation}
where $|B| = \int_B dA$ for any measurable set $B \subset \bbD$, 
and $	D(z,r) = \{ w\in\bbD : \beta(z,w) < r\}$ 
is the disc with center $z$ and radius $r$ in the  Bergman metric 
$\beta: \bbD \times \bbD \to \bbR^+$. Moreover, $\widehat f_r$ is the average function defined by
\bea
	\widehat f_r(z) = \frac{1}{|D(z,r)|} \int\limits_{D(z,r)} f\, dA
	\label{2.1a}
\eea
for $z\in \bbD$. If, in addition,
$$
	\lim_{|z|\to 1} \frac1{|D(z,r)|} \int\limits_{D(z,r)} |f (\zeta) - \widehat f_r(z)|^p \, dA(\zeta)  = 0,
$$
we say that $f$ is in $\VMO^p_r$. These definitions are independent of $r$, and we write simply $\BMO^p$ and $\VMO^p$ for the two spaces when $r=1$. 

To decompose $\BMO^p$ and $\VMO^p$ into smaller spaces, we define the \emph{oscillation} $\omega(f)$ of a continuous function $f$ by
\begin{equation}\label{e:osc}
	\omega(f)(z) = \sup_{w\in D(z,1)} |f(z)-f(w)|
\end{equation}
for $z\in\bbD$. (We fix here the radius of the hyperbolic disc in 
order to keep the notation simple.) The space of \emph{bounded 
oscillation} $\BO$ consists of all continuous functions $f$ for which 
$\omega (f)\in L^\infty$. We say $f\in \BO$ is of \emph{vanishing 
oscillation} and write $f\in \VO$ if $\omega(f)(z) \to 0$ as $|z|\to 1$. 
The spaces $\BA^p$ and $\VA^p$ of functions $f$ of \emph{bounded} or 
\emph{vanishing average} are defined by requiring that $\widehat{|f|^p_1} 
\in L^\infty$ or $\widehat{|f|^p_1}(z)\to 0$ as $|z|\to 1$, respectively. 
These spaces provide the useful decompositions
\begin{equation}\label{2.4}
	\BMO^p = \BO + \BA^p\quad{\rm and}\quad \VMO^p = \VO + \VA^p
\end{equation}
for $p\ge 1$, which can be obtained by setting
$ 
	f = \widehat f_1 + ( f-\widehat f_1 ) .
$ 
For the proofs and further details, see~\cite{MR1178032}.

We will also need the Berezin transform $\widetilde f$ of $f$, which plays an important role in characterizations of various properties of Toeplitz operators. Given $f\in L^1$, the Berezin transform is defined by setting
(see \cite{Z}, Sect. 6.3.)
\begin{equation} \label{2.7}
\widetilde{f}(z) = \frac{\langle fK_z , K_z \rangle_{L^2}
}{\langle K_z , K_z\rangle_{L^2}} = \int\limits_{\bbD} f |k_z|^2 \, dA
= \int\limits_{\bbD} (f \circ \phi_z)
(\zeta) \, dA(\zeta),
\end{equation}
where $k_z  = K_z \Vert K_z \Vert_2^{-1}$ and $\phi_z (w)= \frac{z-w}{1-w\bar z}$ is the M\"obius transform 
interchanging $0$ and $z$. It is a direct consequence of the
definition that the Berezin transform of any function $f \in L^{\infty}$ 
is bounded and continuous.

Recall that a bounded linear operator $T$ on a Banach space $X$ is 
called a \emph{Fredholm} operator, if its kernel ker\,$T$ and cokernel are both finite dimensional. 
If $T$ is Fredholm on $X$, its \emph{index} is defined as the number
$
	\ind T = \dim \ker T - \dim(X/T(X)).
$
The \emph{essential spectrum} ${\rm spec}_{\rm ess}(T)=
{\rm spec}_{\rm ess}(T: X \to X )$ of $T$ is defined by
$$
	{\rm spec}_{\rm ess}(T) 
	= \{\lambda\in\bbC : T-\lambda\ \text{\rm is not Fredholm}\},
$$
which is clearly contained in the spectrum ${\rm spec}(T) = 
{\rm spec}(T:X \to X )$ of $T$. Similarly, the {\it essential norm} is the expression
$ \|T\|_{\rm ess} = \inf\limits_{K} \|T + K\|,$ 
where on the right we have the operator norm of $T+K$ and the infimum is taken 
over all compact operators $K : X \to X$. Clearly, $\|T\|_{\rm ess} \leq \|T\|$.

We still need one more notion to formulate the next theorem. The {\it Stone-\v{C}ech compactification} $\beta \bbD$  of $\bbD$ is defined by its universal property that any 
continuous map $f$ from $\bbD$ to a compact Hausdorff space $K$ can be 
uniquely extended to a continuous map $f: \beta\bbD \to K$. Here, we do 
not distinguish between $f$ and its extension to $\beta\bbD$. Note that 
$\beta\bbD$ can be realized as the maximal ideal space of bounded 
continuous functions defined on $\bbD$. Every maximal ideal corresponds 
to a point in $\beta\bbD$ via evaluation. See e.g. \cite{StroeZhe} for 
the use of $\beta \bbD$ in the topic under consideration.

The following result is known 
and our aim is to extend it  for a  larger symbol class, see Section 
\ref{sec4}.  
\BET 
\label{th2.1}
Let $1 < p < \infty$ and $f \in L^\infty \cap \VMO^1 $.  

\noindent $(i)$ We have
\begin{equation}\label{2.8}
	{\rm spec}_{\rm ess}(T_f : A^p \to A^p ) 
	= \bigcup\limits_{y \in \beta\bbD \setminus \bbD} \widetilde{f}(y)
	= \widetilde{f}(\beta\bbD \setminus \bbD),
\end{equation}
where $\widetilde{f}$ denotes the extension of the Berezin transform of $f$ to the Stone-$\check{\mbox{C}}$ech compactification $\beta\bbD$ of $\bbD$.

\noindent $(ii)$ If $T_f$ is Fredholm on $A^p $, then the index of $T_f$ equals the negative winding of $\widetilde f |_{\{|z|=r\}}$, where $r$ is sufficiently close to $1$.
\ENT

Formula \eqref{2.8} was obtained for the classical Bergman space 
$A^2(\bbD)$ in \cite{Zhu87}. For arbitrary $1 < p < \infty$, \eqref{2.8} was proved 
in~\cite{PeVi} using elementary  methods and in~\cite{Hagger1,Hagger2} 
using techniques with band-dominated operators.
The index formula stated in (ii) was proved 
in~\cite{StroeZhe} for the Hilbert space $A^2 $ and can be treated 
analogously for other values of $p$. 
We remark that, as shown in \cite{StroeZhe}, formula \eqref{2.8} 
and claim (ii) also hold for Toeplitz operators $T_f : A^2 \to A^2$ with 
symbols  in a larger algebra $\cA$ consisting of bounded functions $f$  
such  that $H_f : A^2 \to L^2$ is compact (see Section 5 and
Theorems 19, 24 of the citation). Our generalization is formulated 
in Corollary \ref{cor4}, and it involves also unbounded symbols among 
other things. 

\medskip

\section{Weak $BMO$-type conditions on the unit disc} 
\label{sec3}

In this section we introduce 
the weak $\BMO$-type condition which is interesting in itself and
may be applied to other considerations as well. 
\BED
\label{def0.3}
For all $z = r e^{i \theta} \in \bbD$ with $r \in [0,1)$ and  
$\theta \in [0, 2 \pi)$ we denote
\bea
B(z) =\{ \rho e^{i  \phi} \ | \
r \leq \rho \leq 1 -  \frac12  (1-r)\  ,
\ \theta \leq \phi \leq \theta  +    \pi  (1-r )  \} .
\label{0.4}
\eea
We denote, for  $\zeta = \tilde r  e^{i \tilde \theta } \in B(z)$,
\bea
B(z,\zeta) =\{ \rho e^{i  \phi} \ | \
r \leq \rho \leq \tilde r \  ,
\ \theta \leq \phi \leq \tilde \theta   \}  
\label{0.4a}
\eea
and 
\bea
\widehat f(z,\zeta) := \frac1{|B(z)|} \int\limits_{B(z,\zeta)}  f \,  dA ,
\label{0.6}
\eea
where $f \in L^1_{\rm loc}$. 
\END

\BER
\label{rem3.1}
$(i)$ 
In \eqref{0.4} it may of course happen that $ \theta  +    \pi  (1-r )
> 2 \pi$. This does not harm the definition of the set $B(z)$, but
in \eqref{0.4a} and \eqref{0.6} and in the sequel in all similar places 
we must  agree that the relation $\zeta =\tilde r  e^{i \tilde \theta } 
\in B(z)$ is understood to imply $\tilde \theta \in  [\theta, 
\theta  +    \pi  (1-r )]$, even if  $ \theta  +    \pi  (1-r ) > 2 \pi$.

Using this convention, we define 
for $\zeta_j = \rho_j e^{i \phi_j} \in B(z)$, $j=1,2$, the notion 
$\zeta_1 \precsim \zeta_2$, if $\rho_1 \leq \rho_2$ and $\phi_1
\leq \phi_2$ .

\noindent $(ii)$ In \cite{TV} and \cite{TV2} we specified a 
sequence $(z_n)_{n=1}^\infty \subset \bbD$  such that the 
corresponding sets $B_n := B(z_n)$ form an essentially disjoint 
union of the disc $\bbD$. Here, sets are called essentially disjoint,
if they are disjoint save possibly their boundaries. We will use this 
decomposition later. 
\ENR

In the following we will study symbols $f$ 
for which  there exists a constant $C >0$ such that
\bea
|\widehat f(z,\zeta) |\leq C  \label{4}
\eea
for all $z \in \bbD$ and all $\zeta \in B(z)$. 
By Theorem 2.3 in \cite{TV},  if \eqref{4} holds for the symbol $f$, the 
Toeplitz operator $T_f : A^p \to A^p$, defined as the limit 
\bea
\lim\limits_{\rho \to 1^-} T_{\chi_\rho f}   \label{4re}
\eea 
converging in the strong operator topology,  is bounded. Recall that 
$\chi_\rho $ denotes the characteristic function of the set 
$\bbD_\rho = \{ w \in \bbD \, : \, |w| \leq \rho \}$ with $\rho < 1$. 
Also, according to \cite{TV}, Theorem 2.5, if   $f \in L_{\rm loc}^1$ is 
such that in addition to  \eqref{4} there holds
\bea
\lim_{|z| \to 1 } \frac{1}{|B(z)|} \sup_{\zeta 
\in B(z)} 
\Big| \int\limits_{B(z,\zeta) }f  \, dA \Big|
 = 0,
\label{4a}
\eea
then  $T_f : A^p \to A^p$ is compact. 

Given a precompact subsets $K \subset \bbD $ with $ |K| > 0$,  
we denote the average of $f$ in $K$ by (cf. \eqref{2.1a})
\bea
\widehat f_K = \frac{1}{|K|} \int\limits_K f \,  dA .
\eea
For all $f \in L_{\rm loc}^1$ we also define 
the average function $\widehat f : \bbD \to \bbC$  by
\bea
\widehat f (z) = \widehat f_{B(z)} 
\ , \  \ z  \in {\bbD} .     \label{16}
\eea
Then, we have $\widehat f \in C({\bbD})$ (the space of 
continuous functions on the open disc).

\BER 
\label{rem3.3} 
Note that the more standard definition of the average function
using  a hyperbolic disc instead of $B(z)$ was already
introduced in Section \ref{sec2}, and we keep the difference in the 
notation,  $\widehat f$ vs. $\widehat f_r$, to indicate this. We will need the 
present definition of $\widehat f$ for technical reasons. Moreover,
it is easy to see that in the definitions of the spaces
BO and VO one can replace the sets $D(z,1)$ by the sets
$B(z)$ without changing the concept. This follows from the
simple geometric observation that there exists a number $N \in \bbN$
such that for all $z\in \bbD$, the set $D(z,1)$ is contained
in the union of at most $N$ sets $B(w)$ and conversely, $B(z)$ is
contained in the union of at most $N$ sets $D(w,1)$. For similar reasons, 
in the  definition of the spaces ${\rm BA}^p$ and ${\rm VA}^p$ the 
average functions 
$\widehat{|f|_1^p}$ could be replaced by the average functions
$\widehat{|f|^p}$. We leave the details of the proofs for these claims
to the reader.
\ENR

\BED \label{def4}
Let us consider  functions  $f \in L_{\rm loc}^1$ and define the 
following $\BMO$-type condition 
\bea 
\Vert f \Vert_{\BWMO} := \sup_{z \in \bbD} \frac{1}{|B(z)|} \sup_{\zeta \in B(z)} \Big| \int\limits_{B(z,\zeta) }   
\big( f(\xi)-  \widehat f(z) 
\big) dA (\xi) \Big| < \infty  .  \label{10}
\eea
We refer to this definition as the $\BWMO$-condition (for 
``bounded weak mean oscillation,'' not to be confused with the existing term of ``weak $\BMO$'' in the literature).

We also introduce the following $\VWMO$-condition (for vanishing
weak mean oscillation) for a function $f \in L_{\rm loc}^1$,
\bea
\lim_{|z| \to 1 } \frac{1}{|B(z)|} \sup_{\zeta \in B(z)} \Big| \int\limits_{B(z,\zeta) } 
\big( f (\xi)-  \widehat f(z)  \big) dA(\xi)  \Big| = 0.
\label{14}
\eea
\END

It is easy to see that the expression $\Vert \cdot \Vert_{\BWMO}$
is a seminorm for example in the space of bounded continuous 
functions in $\bbD$ and that $\Vert f \Vert_{\BWMO} = 0$, if
and only if $f$ is constant. To see the latter statement, if $f$ 
is a non-constant, bounded and continuous function in $\bbD$, 
we pick up a point $z \in \bbD$ such that, for a neighborhood $U$ of $z$,
the function $f$ is non-constant in $U \cap B(z)$. Then, it is clear that the 
expression 
$$
\int\limits_{B(z,\zeta) } \big(f -  \widehat f(z)  \big)\, dA 
$$
cannot be constant for $\zeta \in U \cap B(z)$, which 
implies that $\Vert f \Vert_{\BWMO} \not= 0$.

We will use  the following fact. 

\BEL \label{lem5}
Assume \eqref{10} holds for a function $f \in L_{\rm loc}^1 (\bbD)$. 
Let $z \in \bbD$  be arbitrary and assume that  the points
$\tilde z, \zeta$, $\tilde z \precsim \zeta$, belong to $B(z)$ {\rm (}thus, 
$B(\tilde z, \zeta)  \subset  B(z)${\rm )}. Then, we have
\bea
\frac{1}{|B(z)|}  \Big| \int\limits_{B(\tilde z,\zeta) }   
\big( f (\xi) -  \widehat f(z)  \big) dA (\xi)
 \Big| \leq C \Vert f \Vert_{\BWMO}  \label{30}
\eea
for some constant $C>0$. 
\ENL

Proof. We start by the following  elementary geometric
observation: if $g \in L_{\rm loc}^1(\bbD)$, $z = r e^{i \theta} 
\in \bbD$ and  
$\zeta_1, \zeta_2 \in B(z)$  are such that  $z \precsim \zeta_1 
\precsim\zeta_2$ 
then the integral over the set $B(\zeta_1,\zeta_2)$ can be presented  as
\bea
\int\limits_{B (\zeta_1, \zeta_2)} \!\!\!\!\!
g \, dA  = \sum_{j=1}^4 \gamma_j \!\!\!\!\!
\int\limits_{ B(z,w_j ) } \!\!\!\!\! 
 g \,  dA ,
\label{19}
\eea
where $\gamma_j \in \{-1,1 \}$ and  $w_j$, $j = 1, \ldots, 4$ 
are some points in $B(z)$ with $z \precsim w_j \precsim\zeta_2$.  
Indeed, if $\zeta_j = \rho_j e^{i \phi_j}$, $j=1,2$ we choose 
\bea
w_1 = \zeta_2, \ w_2 = \rho_1 e^{i \phi_2}, \ 
w_ 3 = \rho_2 e^{i \phi_1} , \ w_4 = \zeta_1 . 
\eea
Then, we have
\beas
B(\zeta_1, \zeta_2) = B(z,w_1) \setminus
\Big( B(z,w_3 ) \cup\big(  B(z,w_2 ) \setminus  B(z,w_4) \big) \Big)  
\eeas
and formula \eqref{19} follows from this, since the sets 
$  B(z,w_3 )$ and $ B(z,w_2 ) \setminus  B(z,w_4)$ are essentially
disjoint. 

We now apply formula \eqref{19} to the integral in \eqref{30}
and obtain 
\bea
\int\limits_{B(\tilde z,\zeta) }   \!\!\!\!\!
\big( f (\xi) -  \widehat f(z)  \big) dA (\xi)
= \sum_{j=1}^4 \gamma_j 
\int\limits_{B(z,w_j) }   \!\!\!\!\!
\big( f(\xi) -  \widehat f(z)  \big) dA(\xi)  \label{48}
\eea
for some points $w_j \in B(z)$. The bound \eqref{30} follows
from this and the triangle inequality, since  \eqref{10} implies, 
for all $j$, 
\bea
\hskip3cm 
\frac{1}{|B(z)|}
\Big| \int\limits_{B(z,w_j) }   \!\!\!\!\!
\big( f (\xi) -  \widehat f(z)  \big) dA (\xi) \Big| \leq \Vert f \Vert_{\BWMO}.
\hskip3cm  \Box   \label{49}
\eea

\BER \label{rem7}
If, in addition, \eqref{14} holds for $f \in L_{\rm loc}^1$, then 
\bea
\lim\limits_{|z| \to  1} 
\frac{1}{|B(z)|}  \sup\limits_{\stackrel{\scriptstyle \tilde z , 
\zeta \in B(z)}{ \tilde z \precsim \zeta} } 
\Big| \int\limits_{B(\tilde z, \zeta) }   
\big( f (\xi) -  \widehat f(z)  \big) dA(\xi) \Big|
= 0 ,   \label{31}
\eea 
This is so since the assumption \eqref{14} implies that 
in \eqref{48}, \eqref{49} an arbitrarily 
small multiplier $\varepsilon > 0$ can be added to the right hand
side, if $|z| $ is close enough to 1. Following the proof, also the right
hand side of \eqref{49} then can be multiplied with $\varepsilon$, which 
proves the claim. 
\ENR

\BEC
\label{lem3}
There is a  constant $C > 0$ such that, 
if $f \in L_{\rm loc}^1$ satisfies \eqref{10}, then  we have, 
for all $z \in \bbD$ and  $\tilde z, \zeta \in B(z)$  with
$\tilde z \precsim \zeta$,
\bea
\big| \widehat f (z)  - \widehat f_{B(\tilde z, \zeta)}
\big| \leq C \Vert f \Vert_{\BWMO} , \label{25a}  
\eea
provided that the points $z, \tilde z, \zeta$ in addition satisfy 
\bea
 \frac{|B(z) |}{|B(\tilde z, \zeta)|} \leq 2.   \label{25b} 
\eea

If in addition  \eqref{14} holds for $f \in L_{\rm loc}^1$, then we have  
\bea
\lim_{|z| \to 1} 
\sup\limits_{\tilde z , \zeta } 
\big| \widehat f (z)  - \widehat f_{B(\tilde z, \zeta)}
\big| = 0 , \label{25c}  
\eea
where the supremum is taken over all  $\tilde z, \zeta \in B(z)$  with
$\tilde z \precsim \zeta$ satisfying also \eqref{25b}.
\ENC

Proof. To prove the first statement, we fix $z \in \bbD$
and using Lemma \ref{lem5} obtain for the average of 
$f $ over the set $B(\tilde z, \zeta)$
\bea
& & | \widehat f_{B(\tilde z, \zeta)} - \widehat f(z)  | = 
\frac{1}{|B(\tilde z, \zeta)|}
\Big| \!\!\! \int\limits_{B(\tilde z, \zeta)} \!\!\! 
( f (\xi) - \widehat f(z) \, dA (\xi) \Big| 
\rowleq
C \Vert f \Vert_{\BWMO} \frac{|B(z) | }{|B(\tilde z, \zeta)|} \leq 
C' \Vert f \Vert_{\BWMO} , \label{26}
\eea 
where the constant $C$ indeed does not depend on $z, \tilde z, \zeta$.

If \eqref{14} is also satisfied by $f$ and $\varepsilon > 0$ is 
arbitrary, then  we can use Remark \ref{rem7}. This allows a bound 
where $\varepsilon$ is multiplying the right hands sides of \eqref{26}. 
\ \ $\Box$ 

\bigskip

We remind  that $\widehat f_r \in \BO$ or $\widehat f_r \in \VO$ whenever $f\in \BMO^1$ or $f\in \VMO^1$, respectively (see~\eqref{2.4}).  We improve on these in the following theorem, which is the main result of this section. See Remark \ref{rem3.3} for some relevant explanations, and also \eqref{2.4}. 

\BET
\label{prop10}
If $f$ satisfies the $\BWMO$-condition  \eqref{10}, then the function $\widehat f$ {\rm (}see \eqref{16}{\rm )} belongs to $\BO $, and 
if $f$ satisfies the $\VWMO$-condition \eqref{14}, then $\widehat f \in \VO$.
\ENT

Proof. Let us again fix  $z \in \bbD$ and consider 
$w$ such that hyperbolic distance $\beta( z,w) \leq \delta$; 
by choosing $\delta > 0$ small  enough (independently of $z,w$)
it is clear 
that the  intersection $B(z) \cap B(w)$ can be written as $B(\tilde z,\zeta)$ 
and $B(\tilde w, \xi)$ for some points $\tilde z, \zeta \in B(z)$ 
and $\tilde w , \xi \in B(w)$ such that both
\bea
 \frac{|B(z) |}{|B(\tilde z, \zeta)|} \leq 2
 \ \ \mbox{and} \ \ 
  \frac{|B(w) |}{|B(\tilde w, \xi)|} \leq 2  .
\eea
For example, if $z = re^{i \theta}$ and  $w  \in B(z)$, we 
choose $\tilde z = \tilde w = w$ and 
$$
\zeta = \xi = \big(1 - 2^{-1}(1-r) \big)e^{ i (\theta + \pi (1-r))},
$$
see the definition of the set $B(z)$ in \eqref{0.4}. 
We get 
\bea
& & | \widehat f (z) - \widehat f(w)| 
\leq | \widehat f (z) - \widehat f_{B(\tilde z,\zeta)}| + 
|  \widehat f_{B(\tilde w,\xi)} - \widehat f(w)|
\eea
and both of these terms can be bounded by a constant times $\Vert f 
\Vert_{\BWMO}$, by Corollary \ref{lem3}. The second claim follows then 
from the last part of Corollary \ref{lem3}.
\ \ $\Box$

\section{Applications to spectra and Fredholm properties}
\label{sec4}

In order to characterize the essential spectra of Toeplitz operators with $\VWMO$-symbols, we 
recall that the Berezin transform $\widetilde{T}$ of a bounded linear operator $T : A^p \to A^p$ is given by 
\[
\widetilde{T}(z) := \frac{\langle TK_z , K_z \rangle_{L^2}}{\langle K_z , K_z\rangle_{L^2}} , \quad z \in \bbD,
\]
which is well-defined  because $K_z \in A^p$ for all $p \in (1,\infty)$ 
(cf.~\cite{Z}, Ch.~6).  Moreover, $\widetilde{T}$ is always a bounded 
continuous function on $\bbD$. For $f \in L^1$ such that the Toeplitz operator $T_f$ 
is bounded in $A^2$  we get $\widetilde{T_f} = \widetilde{f}$; see
\cite{Z}, p. 165. We will need a generalization.

\BED
\label{lem4.1}
If  $f \in L_{\rm loc}^1(\bbD)$ satisfies \eqref{4}, then we define 
its Berezin transform by $\widetilde f := \widetilde{T_f}$.
\END

Recall that condition \eqref{4} implies the boundedness of $T_f$
in the space $A^2$, hence, the definition coincides with the conventional
one in the case mentioned above. 
We note that the Berezin transform is again independent of $p$. 
In the following,  we still denote the unique extension 
of the Berezin transform of $f$ to $\beta\bbD$ by $\widetilde{f}$. Let $C_{\partial}(\bbD)$ denote the 
set of continuous functions on $\bbD$ that have zero limits at the 
boundary (equivalently: continuous functions on $\beta\bbD$ that vanish 
on $\beta\bbD \setminus \bbD$).

The following corollary includes generalizations of a number of 
known results to new symbol classes. Concrete symbols satisfying 
the assumptions of the corollary are considered in Example \ref{ex5}.  
Estimates for the essential norm of operators in the Toeplitz algebra on 
$A^p$ were previously obtained by Su\'arez~\cite{Suarez} and on weighted 
Bergman spaces $A^p $ by Mitkovski, Su\'arez and Wick~\cite{MiSuWi}, 
which were further improved and extended to bounded symmetric domains 
in~\cite{Hagger2}. We emphasize that the essential norm was previously 
computed exactly only when $p=2$, while it remains an open problem to 
find a sharp constant for the upper estimate for the other values of $p$.
Part (b) of Corollary \ref{cor4} was proved in~\cite{Zorboska} for 
symbols $f$ in  $BMO^1$ that satisfy the condition $\tilde f \in 
L^\infty\cap VO$ (or  equivalently, $\hat f \in L^\infty \cap VO$). 

\BEC
\label{cor4} 
Assume that the symbol $f \in L_{\rm loc}^1(\bbD)$ 
satisfies $\VWMO$-condition \eqref{14}
and that the average function ${\widehat f}$ belongs to 
$L^\infty$. Then, the Toeplitz operator
$T_f: A^p \to A^p$ is bounded for all $1 < p < \infty$ 
and there holds
\begin{itemize}
	\item[(a)] $\widetilde{f} - \widehat{f} \in C_{\partial}(\bbD)$,\vspace{0.2cm}
	\item[(b)] ${\rm spec}_{\rm ess}(T_f) = \widehat{f}(\beta\bbD \setminus \bbD) = \widetilde{f}(\beta\bbD \setminus \bbD)$,\vspace{0.2cm}
	\item[(c)] $\max\limits_{y \in \beta\bbD \setminus \bbD} |\widetilde f(y)| = \max\limits_{y \in \beta\bbD \setminus \bbD} |\widehat f(y)| \leq \Vert T_f \Vert_{\rm ess} \leq \Vert P \Vert \max\limits_{y \in \beta\bbD \setminus \bbD} |\widehat f(y)|,$\vspace{0.2cm}
	\item[(d)] $\Vert T_f \Vert_{\rm ess} = \max\limits_{y \in \beta\bbD \setminus \bbD} |\widehat{f}(y)| = \max\limits_{y \in \beta\bbD \setminus \bbD} |\widetilde{f}(y)|$ for $p = 2$.\vspace{0.2cm}
\end{itemize}
Moreover, $T_f$ is Fredholm if and only if $T_{{\widehat f}}$ is if and only if $T_{{\widetilde f}}$ is.
\ENC

Proof. By  Theorem \ref{prop10} we have ${\widehat f} \in \VO$ 
which clearly  implies that 
\bea
\lim_{|z| \to 1 } \frac{1}{|B(z)|} \sup_{\zeta \in B(z)} \Big| \int\limits_{B(z,\zeta) } 
\big( \widehat f -  \widehat f(z)  \big) dA  \Big| = 0.
\label{140}
\eea
Combining this with the assumption that $f$ satisfies \eqref{14} we see
that the function $f -\widehat f $ satisfies condition \eqref{4a}. This implies that the Toeplitz operator
$ 
T_{f- {\widehat f}} = T_f - T_{{\widehat f}} 
$
is compact while $T_{{\widehat f}}$ is bounded by assumption. It also implies 
that $f$ satisfies \eqref{4}, which means that $\widetilde{f}$ is well-defined, by Definition \ref{lem4.1}. Further, by \cite{Zhu87},
we have $\widehat{f} - \widetilde{\widehat{f}\ } \in C_{\partial}(\bbD)$ and
\[
\widetilde{f} - \widehat{f} = \widetilde{\widehat{f}\ } - \widehat{f} + (f - \widehat{f})^{\sim} \in C_{\partial}(\bbD),
\]
where $(f - \widehat{f})^{\sim} \in C_{\partial}(\bbD)$ follows from \cite[Theorem 9.5]{Suarez}. This proves (a). 

Also, the case (b) follows from the above observations and Theorem \ref{th2.1}.
Moreover, we find that it is enough to prove (c) and (d) for the 
function $\widehat f$ instead of $f$. Hence, by Theorem \ref{prop10}
and a redefinition of the notation, we may assume that
$f \in L^\infty \cap \VMO^1$ for the rest of the proof. 
First, for all $z \in \bbD$ we denote by $U_z : L^p \to L^p$ the 
surjective isometry (reflection)
\begin{equation}\label{U_z}
	(U_zf)(w) = f(\phi_z(w))\frac{(1-|z|^2)^{2/p}}{(1- 
	z \bar w)^{4/p}}, \ \ \ w\in \bbD.
\end{equation}
Notice that $U_z^{-1} = U_z$ and $U_zM_fU_z^{-1} = M_{f \circ \phi_z}$, where $M_f$ is the multiplication operator defined by $M_fg = fg$.

Let $z \in \beta\bbD$ and choose a net $(z_{\gamma})$ in $\bbD$ 
that converges to $z$. We note that the operator
$U_{z_{\gamma}}T_fU_{z_{\gamma}}^{-1}$ converges strongly to $T_{h_z}$,
where $h_z$ is a bounded and analytic function on $\bbD$. This follows 
from Lemma~2.1 and Proposition~2.2 of \cite{HaVi}. 
By  Theorem 22 and formula (4.1) of \cite{Hagger2}, we have
\[
\sup\limits_{z \in \beta\bbD \setminus \bbD} \Vert T_{h_z} \Vert 
\leq \Vert T_f \Vert_{\rm ess} \leq \Vert P 
\Vert
\sup\limits_{z \in \beta\bbD \setminus \bbD} \Vert T_{h_z} \Vert .
\]
As $h_z$ is bounded and analytic, we have $\Vert T_{h_z} \Vert \leq 
\Vert h_z \Vert_{\infty}$. Moreover, it is well-known that
${\rm spec} (T_{h_z}) = {\rm clos}(h_z(\bbD))$ (``clos'' denotes 
the closure of the set). 
Indeed, $T_{(h_z-\lambda)^{-1}}$ is an inverse of $T_{h_z - \lambda}$ 
if $\lambda \notin {\rm clos}(h_z(\bbD))$ and conversely
$\overline{h_z(w)}$ is an eigenvalue of $T_{h_x}^*$ for every 
$w \in \bbD$ since 
\[
T_{h_z}^*K_w = P (\overline{h_z}K_w) = \overline{h_z(w)}K_w .
\]
Since $\Vert T_{h_z} \Vert \geq w$ for all $w \in {\rm spec} (T_{h_z})$,
we  get $\Vert T_{h_z} \Vert  = \Vert h_z \Vert_{\infty}$.
Moreover, 
\begin{equation*} 
\lim\limits_{z_{\gamma} \to x} \tilde{f}(\phi_{z_{\gamma}}(w)) = \lim\limits_{z_{\gamma} \to x} \widetilde{f \circ \phi_{z_{\gamma}}}(w) = \widetilde{h_x}(w) = h_x(w)
\end{equation*} 
implies
\[
\sup\limits_{x \in \beta\bbD \setminus \bbD} 
\Vert h_z \Vert_{\infty} = \max\limits_{y \in \beta\bbD \setminus \bbD} |\tilde{f}(y)|.
\]
Combining these estimates, we obtain (c) and (d).  \ \ $\Box$

\begin{corollary}
If $f\in L^1(\bbD)$ satisfies the $\VWMO$-condition 
\eqref{14}, $\widehat{f}$ is bounded and $T_f$ is Fredholm on $A^p$ for
$1 < p < \infty$, then 
$$
	\ind T_f = \ind T_{{\widehat f}} = 	-\wind \widehat f |_{\{|z|=r\}}
$$
where $\wind$ denotes the winding number and $r$ is sufficiently close to $1$. The same statement also holds if $\widehat{f}$ is replaced by $\widetilde{f}$.
\end{corollary}

Proof. Notice first that 
$$
	\ind T_f = \ind T_{{\widehat f}} + \ind T_{f-{\widehat f}} = \ind T_{{\widehat f}}
$$
because $T_{f-{\widehat f}}$ is compact. By 
Theorem \ref{prop10}, ${\widehat f}\in \VO$, and the rest now follows from Theorem~\ref{th2.1} (ii) and Corollary \ref{cor4}. $\Box$

\bigskip

In the next result we need to assume $f \in L^1(\bbD)$ instead of 
$f \in L_ {\rm loc}^1(\bbD)$, although some  implications would
hold even for locally integrable symbols. 

\BEC
\label{cor5} 
Assume that the symbol $f \in L^1(\bbD)$ 
satisfies $\BWMO$-condition \eqref{10}. Then, the following are
equivalent:

\smallskip

\noindent $(i)$  $T_f : A^p \to A^p$ is bounded for some $1 < p < 
\infty$,

\smallskip

\noindent $(ii)$  $T_f : A^p \to A^p$ is bounded for all $1 < p < 
\infty$,

\smallskip

\noindent $(iii)$ $\widehat{f} \in L^\infty(\bbD)$, 

\smallskip

\noindent $(iv)$  $\widetilde{f} \in  L^\infty(\bbD)$. 
\ENC

Proof. $(iii) \Rightarrow (ii) \Rightarrow (iv)$: Assume $(iii)$ holds. 
We have, for some  constant $C> 0$ and all $z \in \bbD$ and  $\zeta \in B(z)$, 
\bea
& & \frac{1}{|B(z)|} \Big| \int\limits_{B(z,\zeta)} ( f-  \widehat{f} ) dA 
\Big| 
\leq \frac{1}{|B(z)|} \Big| 
\int\limits_{B(z,\zeta)} \big( f (\xi)
- \widehat{f}(z) \big) dA(\xi) \Big|
\rowpl
\frac{1}{|B(z)|} \Big|\int\limits_{B(z,\zeta)} \big( \widehat{f} (\xi)
- \widehat{f}(z) \big) dA(\xi) \Big|
\leq C , \label{141}
\eea
because the first integral on the right is bounded due to 
$f \in \BWMO$  and the second one due to $\widehat f \in L^\infty$. 
Hence, the function 
$f - \widehat{f}$ satisfies condition \eqref{4}, and thus the Toeplitz
operator $T_{f -\widehat f}$ is bounded in any $A^p$. 
The boundedness of the average function  $\widehat f$ also implies that
the operator $T_{f} =  T_{f - \widehat f} + T_{\widehat f}$ is bounded in 
any $A^p$. The boundedness of $T_f$ yields the boundedness of the Berezin
transform $\widetilde f$, see the beginning of Section \ref{sec4}. 

$(i) \Rightarrow (iii)$:
Assume $p$ is such that $T_f$ is bounded in $A^p$. 
We can use \eqref{141} to see that the operator $T_{f -\widehat f}$ is bounded 
in $A^p$. The boundedness  of the second integral on the right of \eqref{141} 
follows from  $\widehat f \in \BO$, see Theorem \ref{prop10}. 
Next, $T_f$  and $T_{f -\widehat f}$ bounded implies $T_{\widehat f}$ 
bounded, hence the Berezin transform $\widetilde{\widehat f \ }$ is
a bounded function, and this fact does not depend on $p$. 
Since $\widehat f \in \BO \subset \BMO^1$
we obtain that also $\widehat{\widehat f \ } \in L^\infty(\bbD)$ (see Corollary~2.3.(d) of~\cite{Z3}). Also, $\widehat f \in \BO$ implies
$ \big(   f - {\widehat f}\big)^\wedge \in L^\infty(\bbD)$.
We obtain
\beas
\widehat f = \big(   f - {\widehat f}\big)^\wedge
+   \widehat{\widehat f \ } 
\in L^\infty(\bbD). 
\eeas

$(iv) \Rightarrow (iii)$
It suffices to show that $\widetilde f \in L^\infty$ implies 
$\widetilde{\widehat f \ } \in L^\infty$. Once we have this,
we deduce as in the previous item that   $\widehat{\widehat f \ } 
\in L^\infty(\bbD)$ and then $ \big(   f - {\widehat f}\big)^\wedge \in 
L^\infty(\bbD)$ and finally $\widehat f  \in L^\infty(\bbD)$.

We utilize the decomposition of the disc into the hyperbolic
sets $B_n = B(z_n)$, mentioned in Remark \ref{rem3.1}.$(ii)$ and
write
\bea
& & \widetilde{\widehat f \ } (z) = \int\limits_\bbD 
\widehat f  
|k_z|^2  dA
= 
\sum_{n \in \bbN} \int\limits_{B_n}
\widehat f(w) |k_z(w)|^2 dA(w).
\label{20a}
\eea
We note that since $\widehat f \in \BO$, the modulus of the
expression 
\beas
\widehat f (w) - \widehat f(z_n)
=: F_n(w)
\eeas
is bounded  for every $w \in B_n $ by a constant $C_1 >0$ independent of 
$w$ or $n$. Hence, \eqref{20a} equals
\bea
& & \sum_{n \in \bbN} \int\limits_{B_n}
\widehat f (z_n) |k_z(w)|^2 dA(w)
+ \sum_{n \in \bbN} \int\limits_{B_n}
F_n |k_z|^2 dA ,  \label{20c}
\eea
where 
\beas
\Big| \sum_{n \in \bbN} \int\limits_{B_n}
F_n |k_z|^2 dA \Big| \leq 
\sum_{n \in \bbN} \int\limits_{B_n} C_1 |k_z|^2 dA = 
C_1 \Vert k_z \Vert_ 2^2 = C_1 
\eeas
since $k_z$ is the normalized kernel, see \eqref{2.7}. 
To estimate the first  term  in \eqref{20c} we write 
\bea
G_z(w,\zeta) = |k_z(w)|^2  - |k_z(\zeta)|^2 ,
\eea
which again is bounded by a constant $C_2>0$ independent of $z
\in \bbD$ and  $n\in \bbN$ and $w, \zeta \in  B_n$. We get
\bea
& & 
 \sum_{n \in \bbN} \int\limits_{B_n}
\widehat f (z_n) |k_z(w)|^2 dA(w) 
= \sum_{n \in \bbN} \int\limits_{B_n}
\frac{1}{|B_n|} \int\limits_{B_n} f(\zeta) |k_z(\zeta)|^2 
dA(\zeta) dA(w)
\rowpl
\sum_{n \in \bbN} \int\limits_{B_n} \frac{1}{|B_n|} 
\int\limits_{B_n} f(\zeta)  G_z(w, \zeta) dA(\zeta) dA(w).
\eea
Here, the first term equals $\widetilde f$. The second one is bounded by
\beas
\sum_{n \in \bbN} \int\limits_{B_n} \frac{1}{|B_n|} 
\int\limits_{B_n} C_2 |f(\zeta)| dA(\zeta) dA(w)
\leq 
 C_2 \sum_{n \in \bbN} 
\int\limits_{B_n} |f(\zeta)| dA(\zeta)  \leq 
C_2 \Vert f \Vert_1.
\eeas
This completes the proof of the corollary. \ \ $\Box$

\BEC
\label{cor8}
$(i)$ Assume that the symbol $f \in L^1(\bbD)$  satisfies $\BWMO$-condition 
\eqref{10}. Then, the Toeplitz operator $T_{f - \widehat f}$ is bounded in $A^p$ for all $1 < p < \infty$.

\noindent $(ii)$ If  $f \in L^1(\bbD)$ satisfies $\VWMO$-condition 
\eqref{14}, then $T_{f - \widehat f}$ is compact in $A^p$ for all $1 < p < 
\infty$. If, in  addition,  $T_f$ is  bounded in $A^p$ for some $1 < p < 
\infty$, then also  $T_{f - \widetilde f}$ is compact in all spaces
$A^p$. 

\ENC

Proof. The boundedness of $T_{f - \widehat f}$ follows from the 
proof of Corollary \ref{cor5}, $(i) \Rightarrow (iii)$, and the compactness
from the beginning of the proof of Corollary \ref{cor4}. 
If in addition $T_f$ is bounded, 
Corollary \ref{cor5} yields that  $\widehat f \in L^\infty$. We obtain from 
Corollary \ref{cor4} that the bounded, continuous function $\widetilde f - 
\widehat f $ belongs to $C_\partial (\bbD)$. Hence, the operator
$T_{\widetilde f - \widehat f}$ and thus also $T_{f -\widetilde f}$
are  compact. \ \ $\Box$

\bigskip

\BER \label{rem12}
$(i)$ The proof of Theorem 3.7 in \cite{Zorboska} shows that if $f \in 
\BMO^2$ and the Berezin transform $\widetilde f \in L^\infty \cap 
\VO$, then the Toeplitz operator $T_{f - \widetilde f}$ is compact on 
$A^p$.  
Since in this case $\widetilde f$ is a bounded continuous function,
the results concerning the Berezin transform in (b)--(d) of Corollary \ref{cor4} hold true for such symbols. 

\noindent $(ii)$ It is known that for $f \in \BMO^1$, the compactness
of $T_f: A^p \to A^p$ for some $1 < p< \infty$ is
equivalent to the vanishing of the  Berezin
transform at the boundary, i.e.~$\widetilde f 
\in C_\partial(\bbD)$; see \cite{Z3}, Theorem 3.1. We do not know if 
this result can be generalized for symbols in the BWMO-class.
\ENR

We give one simple application of Theorem~\ref{prop10} to the study of block 
Toeplitz operators. For a suitable $f = (f_{jk})_{j,k=1}^N$ with $f_{jk} 
\in L^1$, $N \geq 2$, denote the block Toeplitz operator  by $T_f : A^p_N
\to A^p_N$, where $A^p_N = \{ g=(g_1, \ldots, g_N) : g_k \in A^p\}$ with 
$\|g\|_{A^p_N} = \max \|g_k\|_p$. More precisely, if 
$$
	g = (g_1,\ldots, g_N) = (g_k)_{k=1,\ldots,N} \in A^p_N,
$$
then
$$
	T_f g = P(fg^T)= P\left( \left(\sum_{k=1}^N f_{jk}g_k\right)_{j=1,\ldots, N}\right) = \left( \sum_{k=1}^N T_{f_{jk}} g_k\right)_{j=1,\ldots, N},
$$
where $g^T$ is the transpose of $g$. Similarly, the Berezin transform of $T_f$ is a matrix operator which is defined componentwise by using the scalar definition in the beginning of Section \ref{sec4}. For further details, see \cite{PeVi}

For scalar symbols  $f,g\in L^\infty$, it is easy to see that
\begin{equation}\label{e:prod}
	T_f T_g = T_{fg} - PM_fH_g.
\end{equation}
Note that the Hankel operator $H_g$ is compact in any $L^p$, if $g \in \VO$, 
see e.g. \cite{HaVi}. 

\begin{proposition}
Let $f = (f_{jk})$ with $f_{jk} \in L^1\,\cap\, \VWMO$ and suppose that 
$\widehat f_{jk} \in L^\infty$. Then $T_f$ is Fredholm on $A^p_N$ if and only if $\det \widetilde f$ is bounded away from zero near $\partial\bbD$.
\end{proposition}

Proof. First, we note that all operators $T_{jk}:A^p \to A^p$ and
$T_f:  A^p_N \to A^p_N$ are  bounded as a consequence of Corollary \ref{cor5}.
Hence, the Berezin transform $\widetilde f = (\widetilde f_{jk })$ is 
a well defined, bounded and continuous matrix function $\bbD \to 
\bbC^{N \times N}$. Moreover, $T_{f-\widetilde f}$ is compact on $A^p_N$
by Corollary \ref{cor8}.

We reduce the matrix-valued case to the scalar case using the following well-
known theorem: if the entries $A_{jk}$ of a bounded linear matrix operator $A$ 
on a product Banach space $X^N$ commute modulo compact operators, then $A$ is 
Fredholm on $X^N$ if and only if $\det A$ is Fredholm on $X$ (see, 
e.g.,~\cite{PeVi}).

Since
$$
	T_f = T_{\widetilde f} + T_{f-\widetilde f}
$$
and $T_{f-\widetilde f}$ is compact on $A^p_N$, it follows that $T_f$ is 
Fredholm if and only if $T_{\widetilde f}$ is Fredholm. By \eqref{e:prod}, if 
the scalar symbols $g,h$ belong to $L^\infty\cap \VO$, then $T_g$ and $T_h$ 
commute modulo compact operators. Now, by Theorem~\ref{prop10}, all 
$\widehat f_{jk}$ belong to $L^\infty\cap \VO$, and the same is true also
for all $\widetilde f_{jk}$, by Corollary \ref{cor4}.(a). We conclude 
by the above mentioned  theorem that $T_{\widetilde f}$, equivalently, $T_f$, 
are Fredholm if and only if $\det T_{\widetilde f}$ is Fredholm. Notice that 
$$
	\det T_{\widetilde f} = \sum_{\sigma\in S_N} \left({\rm sgn}(\sigma)\prod_{j=1}^N T_{\widetilde f_{j,\sigma_j}}\right) = T_{\det \widetilde f} + K,
$$
where $S_N$ is the permutation group and  $K$ is some compact operator. Therefore, because $\widetilde f_{jk} \in L^\infty\cap \VO$ so that $\det \widetilde f\in L^\infty\cap \VO$, the scalar case (see Corollary~\ref{cor4} or \cite{Hagger1, PeVi}) implies that $T_f$ is Fredholm if and only if $\det \widetilde f$ is bounded away from zero near $\partial \bbD$.
\ \ $\Box$

\BEE
\label{ex5}
We present examples of symbols which $(i)$ satisfy the 
$\VWMO$-condition \eqref{14} but which are not in $\BMO^1$, or $(ii)$ bounded 
$\VWMO$-symbols which are not in $\VMO^1$. In the following examples the 
average function $\widehat f$ belongs to $C_\partial(\bbD)$ 
and the operator $T_{\widehat f}$ is compact. Although the example resembles 
that in Remark 2.4  of \cite{TV}, 
the proof is completely different, the reason being that the standard 
definition of $\BMO^1$ involves hyperbolic discs $D(z,r)$ instead of the sets
$B(z)$, and the present technique is more convenient here. We define for all $b \geq \beta >  0$ the function
\begin{equation}
f(r e^{i \theta} ) := \left\{
\begin{array}{ll}
{\displaystyle \frac1{r(1-r)^{b - \beta}}  \sin \frac1{(1-r)^b} } \ , &  \ \ r\geq \frac 12 \\
1 \ , & \ \ r <  \frac 12 .
\end{array} \right.
\end{equation}
Given $z = r e^{i \theta}$ we 
calculate using integration by parts, for all $\zeta = \tilde r e^{i \tilde 
\theta} \in B(z)$, $\xi  = \rho e^{i \phi} \in B(z,\zeta )$, 
\bea
& &  \int\limits_{B(z,\zeta) }
f(\xi) \, dA (\xi)    = \int\limits_\theta^{\tilde \theta}  
\int\limits_{r}^{\tilde r}
(1- \rho)^{1 + \beta} \frac1{(1-\rho)^{b+1} } \sin 
\frac1{(1-\rho)^b}  d \rho d \phi 
\roweq 
- \int\limits_\theta^{\tilde \theta}  
\bigg( \Big[ \frac{(1- \rho)^{1 + \beta} }b \cos \frac1{(1-\rho)^b 
}\Big]_{\rho  = r}^{ \rho= \tilde r} 
\rowpl
\frac{1+\beta}b
\int\limits_{r}^{\tilde r} (1- \rho)^{\beta} \cos \frac1{(1-\rho)^b } d \rho
\bigg) d \phi 
\label{35}
\eea
We estimate $| \cos ( \ldots )| \leq 1$ and take into account the
lengths of the integration intervals, hence, the modulus of
the replacement term is bounded by 
\beas
C  \int\limits_\theta^{\tilde \theta}  (1- r)^{1 + \beta } d \phi
\leq C' (1-r)^{2+ \beta} 
\eeas
and the last integral in \eqref{35} has the same bound. Since $|B(z)| \geq C 
(1-r)^2$, we obtain the following two conclusions: first, $f $ satisfies 
condition \eqref{4a} and second, $|\widehat f| \in C_\partial (\bbD)$. 
These two imply $f \in \VWMO$. 

To see that $f \notin \BMO^1$ we first show that also the average
function $\widehat f_1$, see \eqref{2.1a}, belongs to $C_\partial (\bbD)$. 
To this end we fix $z = r e^{i \theta}$ and consider a set $D(z,1)$ instead of 
$B(z)$. It follows from the definition of the hyperbolic geometry, see 
\cite{Z}, Ch. 6, that for some constant $C> 0$ we have for all $\zeta \in 
D(z,1)$ 
\bea
\frac1C (1 -r ) \leq 1 - |\zeta| \leq C( 1- r)  \ \  \mbox{and} \ 
\sup_{w_1, w_2 \in D(z,1)} | w_1 - w_2| \leq C(1-r)  \label{36}
\eea
It is also obvious that the set $D(z,1)$ can be presented using polar
coordinates as 
\beas
D(z,1)& = & \big\{ \xi = \rho e^{i \phi} \, : \,  \theta - \theta_0 < \phi < 
\theta + \theta_0,
\ r_1 (\phi) < \rho < r_2 (\phi) \big\}
\eeas
for some number $\theta_0> 0$ and functions $r_j : (\theta - \theta_0 , 
\theta + \theta_0) \to (0,1)$. The points  $r_j (\phi) e^{i \phi}$, $j =1,2$, form the 
boundary of the disc  $D(z,1)$; moreover, by \eqref{36}, $\theta_0$ and 
$r_j$ are bounded by $C(1-r)$ and also bounded from below by 
$C'(1-r)$. 

We obtain in the same way as in \eqref {35}
\beas
& &  \int\limits_{D(z,1) } f(\xi) \, dA (\xi)    =
\int\limits_{\theta - \theta_0}^{\theta + \theta_0}  
\bigg( \Big[ \frac{(1- \rho)^{1 + \beta} }b \cos \frac1{(1-\rho)^b 
}\Big]_{\rho  = r_1(\phi)}^{ \rho= r_2(\phi)} 
\rowpl
\frac{1+\beta}b
\int\limits_{r_1(\phi)}^{r_2(\phi)} (1- \rho)^{\beta} \cos \frac1{(1-\rho)^b } d \rho
\bigg) d \phi  ,   
\eeas
and as above we see that the modulus of this is bounded by 
$C(1-r)^{2 + \beta}$. 
Hence, $|\widehat f_1(z)| \leq C(1-r)^\beta$ for all $z = r e^{i \theta}
\in \bbD$ and in particular $\widehat f_1 \in C_\partial (\bbD)$.

It is quite obvious that there is a constant $\delta > 0$ 
such that we have the lower bound 
$|\sin ( 1/(1-r)) |\geq \delta$ for $re^{i\theta}$ in a subset $D_3$ of $D(z,1)$ with measure
at least $|D(z,1)|/ 2$ , thus
$|f(z)| \geq C(1-r)^{\beta-b}$ for $z=re^{i\theta}\in D_3$. As a consequence
\beas
& & \frac{1}{ |D(z,1)|} \int\limits_{D(z,1)} | f(\xi) - \widehat f_1(z)| dA(\xi)
\rowgeq \frac{1}{ |D(z,1)|}\int\limits_{D_3} | f| dA
- \frac{1}{ |D(z,1)|} \int\limits_{D_3} | \widehat f_1(z)| dA (\xi)
\rowgeq 
\frac{C}{ |D(z,1)|} \int\limits_{D_3} \frac{1}{(1-r)^{b - \beta}} dA 
- \frac{C'}{ |D(z,1)|} \int\limits_{D_3} (1-r)^{\beta} dA 
\geq C''(1-r)^{\beta -b} . 
\eeas
In view of the definition of the norm of $\BMO^1$, see \eqref{2.1}, 
we get an example of type 
$(i)$ by taking any parameters $b,\beta$ such that
$b > \beta >0$. Moreover, if $b - \beta < 1$, then there holds $f \in L^1$,  
but if $b - \beta \geq 1$, we only have $f \in L_{\rm loc}^1$. 
To obtain an example of type $(ii)$ one chooses $b = \beta > 0$.  
\ENE

\bigskip

\noindent\textbf{Acknowledgments.} The authors thank Raffael Hagger for useful discussions.


\begin{thebibliography}{99}

\bibitem{Hagger1} R.~Hagger, The essential spectrum of Toeplitz operators on the unit ball. Integral Equations Operator Theory 89 (2017), no. 4, 519--556.

\bibitem{Hagger2} R.~Hagger, Limit operators, compactness and essential spectra on bounded symmetric domains. J. Math. Anal. Appl. 470 (2019), no. 1, 470--499.

\bibitem{HaVi} R. Hagger and J. A. Virtanen, Compact Hankel operators with bounded symbols, J. Operator Theory (in press).


\bibitem{MiSuWi} M.~Mitkovski, D.~Su\'arez and B.~Wick, The essential norm of operators on {$A^p_\alpha(\bbB_n)$}. Integral Equations Operator Theory 75 (2013), no. 2, 197--233.

\bibitem{PeVi} A. Per\"al\"a and J. A. Virtanen, A note on the Fredholm properties of Toeplitz operators on weighted Bergman spaces with matrix-valued symbols. Oper. Matrices 5 (2011), no. 1, 97--106.

\bibitem{StroeZhe} K.~Stroethoff and D.~Zheng, Toeplitz and Hankel operators on Bergman spaces. Trans. Amer. Math. Soc. 329 (1992), no. 2, 773--794.

\bibitem{Suarez} D.~Su\'arez, The essential norm of operators in the {T}oeplitz algebra on {$A^p(\bbB_n)$}. Indiana Univ. Math. J. 56 (2007), no. 5, 2185--2232.

\bibitem{TV} J. Taskinen and J. A. Virtanen, Toeplitz operators on Bergman spaces with locally integrable symbols. Rev. Mat. Iberoam. 26 (2010), no. 2, 693--706.

\bibitem{TV2} J. Taskinen and J. A. Virtanen, On generalized Toeplitz and little Hankel operators on Bergman spaces. Arch. Math. (Basel) 110 (2018), no. 2, 155--166.

\bibitem{YZ} F. Yan and D. Zheng, Bounded Toeplitz operators on Bergman space. Banach J. Math. Anal. 13 (2019), no. 2, 386--406.

\bibitem{Zhu87} K. Zhu, VMO, ESV, and Toeplitz operators on the Bergman space, Trans. Amer. Math. Soc. 302 (1987), no. 2, 617--646.

\bibitem{MR1178032}  K. Zhu, BMO and Hankel operators on Bergman spaces. Pacific J. Math. 155 (1992), no. 2, 377--395.

\bibitem{Z3} N. Zorboska, Toeplitz operators with BMO symbols 
and the Berezin transform, Int. J. Math. Sci. 46 (2003), 
2929--2945.

\bibitem{Zorboska} N. Zorboska, Closed range type properties of Toeplitz operators on the Bergman space and the Berezin transform, Complex Anal. Oper. Theory 13 (2019), no. 8, 4027--4044. 

\bibitem{Z}  K. Zhu, Operator theory in function spaces. Second edition. Mathematical Surveys and Monographs, 138. American Mathematical Society, Providence, RI, 2007.

\end{thebibliography}
\end{document}